\documentclass[12pt]{amsart}
\pdfoutput=1 
\usepackage{etex}
\usepackage{amsmath, amsthm, graphicx, amssymb,verbatim,pst-all,color}
\usepackage[margin=1.0in,footskip=0.45in]{geometry}
\usepackage{mathptmx}
\usepackage[latin1]{inputenc}
\usepackage{multirow}
\usepackage{subfig}
\usepackage{enumerate}
\usepackage{array,nicefrac}
\usepackage{url}
\usepackage{times}
\usepackage{textcomp}
\usepackage{verbatim}

\usepackage{epsf}
\usepackage{stmaryrd}
\usepackage{cases} 
\usepackage{relsize}
\usepackage{mathtools}
\usepackage{amscd, amsfonts}
\usepackage{enumitem}
\usepackage[all]{xy}
\usepackage{lscape}
\usepackage{tikz-cd}
\usepackage{etex,soul}
\usepackage{pgf}
\usepackage{hyperref}
\usepackage{enumitem}
\ifx\JPicScale\undefined\def\JPicScale{1.0}\fi
\unitlength \JPicScale mm
{
\theoremstyle{plain}
   \newtheorem{theorem}{Theorem}[section]
   \newtheorem{proposition}[theorem]{Proposition}     
   \newtheorem{lemma}[theorem]{Lemma}

   \newtheorem{corollary}[theorem]{Corollary}
   
}
{
\theoremstyle{definition}

   \newtheorem{example}[theorem]{Example}
   \newtheorem{definition}[theorem]{Definition}
   \newtheorem{remark}[theorem]{Remark}

}

\newcommand{\CC}{\mathbb{C}}

\newcommand{\ZZ}{\mathbb{Z}}

\newcommand{\QQ}{\mathbb{Q}}

\newcommand{\PP}{\mathbb{P}}

\newcommand{\iso}{\cong}

\newcommand{\wmod}{X_\mathcal{A}[n]}

\newcommand{\bP}{\mathbb{P}}

\author{Patricio Gallardo, Jos\'e Luis Gonz\'alez \and Evangelos Routis}
\address{
{\small Department of Mathematics,
University of California, Riverside,
900 University Ave.
Riverside, CA 92521
Skye Hall}
}
\email{pgallard@ucr.edu}

\address{
{\small Department of Mathematics,
University of California, Riverside,
900 University Ave.
Riverside, CA 92521
Skye Hall}
}
\email{jose.gonzalez@ucr.edu}

\address{
{\small 
Mathematics Institute,
Zeeman Building,
University of Warwick,
Coventry CV4 7AL}
}
\email{evangelos.routis@warwick.ac.uk}

\begin{document}
\title{The Fulton-MacPherson compactification is not a Mori~dream~space}
\maketitle


\begin{abstract}
We show that the Fulton-MacPherson  compactification of the configuration space of $n$ distinct labeled points in certain varieties of arbitrary dimension $d$, including projective space, is not a Mori dream space for $n$ larger than $d + 8$.  

\end{abstract}


\section{Introduction}
Understanding the birational geometry of compact moduli spaces is a fundamental problem within algebraic geometry. A central question is to determine if a moduli space is a Mori dream space, 
that is, a normal, projective, $\mathbb{Q}$-factorial variety with a finitely generated class group and a finitely generated Cox ring.  
This condition implies that the main conjectures and theorems of Mori's program hold for the variety, see  \cite[{1.10-1.11, 2.9-2.14}]{Hu-Keel}.
The above question can be very hard to answer in general. A.-M.~Castravet and J.~Tevelev only recently determined that the moduli space of stable rational pointed curves $\overline{M}_{0,n}$ fails to be a Mori dream space for $n \geq 134$, see \cite{castravet2015overline}. Subsequent works by different teams, including one of the authors, have improved the above lower bound to $n \geq 10$, see \cite{gonzalez2016some}, \cite{hausen2018blowing}. We extend the above results to spaces parametrizing certain pointed higher dimensional objects, which we will describe next. 

In 1994, W. Fulton and R. MacPherson \cite{Fulton-MacPherson} constructed a compactification $X[n]$ of the configuration space of $n$ distinct labeled points in an arbitrary smooth variety $X$, which enjoys several desirable properties. To list a few, $X[n]$ is smooth with normal crossings boundary, it has an explicit blowup construction and its geometric points can be given a tree-like description resembling the one of $\overline M_{0,n}$. 

The above compactification 
is a central object in algebraic geometry, yet its birational geometry is still far from being understood. 
A. Massarenti has studied the automorphism group of $X[n]$, see \cite{massarenti2017biregular}. 
When $X=\PP^1$ the space $X[n]$ can be identified with Kontsevich's moduli space $\overline M_{0,n}(\PP^1,1)$ of stable maps of degree 1 to $\PP^1$ and therefore, recent work of C. Fontanari \cite{fontanari2020overline} shows that $\mathbb{P}^1[n]$ is a Mori dream space whenever $\overline{M}_{0,n+3}$ is. 

In this work, we show that $X[n]$ is not a Mori dream space for a broad range of $X$ over an algebraically closed field of characteristic 0, including $X=\PP^d$ for any $d\geq 1$ and $n\geq d+9$, therefore answering a question asked by Fulton to one of the authors.
\begin{theorem}\label{thm.FM}
 Let $f:X\dashrightarrow \PP^d$ be a birational map where $X$ is a smooth projective variety. Suppose that there exists an open subset $U$ of $X$ and an open subset $V$ of $\mathbb{P}^d$ with $codim(\PP^d\setminus V,\mathbb{P}^d)\geq 2$ such that $f_{|U}$ is an isomorphism onto $V$. Then $X[n]$ is not a Mori dream space for $n\geq d+9$.
\end{theorem}

A major role in the proof of the above theorem is played by certain compactifications $\overline{P}_{d,n}$ of the moduli space of $n$ distinct labeled points in $\PP^d$, up to the usual action of $SL_{d+1}$, introduced by two of the authors in \cite{GR}. These are constructed as GIT quotients of Fulton-MacPherson compactifications for certain $SL_{d+1}$-linearized line bundles $\tilde L_{d,n}$,
\begin{align*}
\overline{P}_{d,n} := 
\mathbb{P}^d[n] / \!\!/_{\tilde L_{d,n}} SL_{d+1}    
&&
\left(
\text{we have }
\overline{P}_{1,n} :=  
\mathbb{P}^1[n] / \!\!/_{\tilde L_{1,n}} SL_{2}(\CC) 
\cong
\overline{M}_{0,n} 
\,\,\text{by \cite[Theorem 3.4]{Hu-Keel}}\right)
\end{align*}
and can be viewed as higher dimensional generalizations of $\overline{M}_{0,n}$, see Section  \ref{subsection.wonderful}.

 The study of toric varieties whose blowup at a general point is not a Mori dream space has attracted a lot of interest recently (see \cite{gonzalez2016some}, \cite{hausen2018blowing}, \cite{gonzalez2019examples}, \cite{he2019mori}, \cite{kurano2019infinitely}, \cite{gonzalez2020curves}) and often plays a pivotal role in determining whether a given variety is a Mori dream space (e.g. \cite{castravet2015overline}, \cite{gonzalez2016some}, \cite{hausen2018blowing}). Equally interesting seems to be the more delicate question of whether such a blowup has a polyhedral pseudoeffective cone (see e.g. \cite{CastravetLafaceTevelevUgaglia2020}, \cite{laface2021intrinsic}). 
 Our Proposition~\ref {proposition.projection} and  Corollary~\ref{corollary.projection} provide a general method for finding such toric varieties via a series of certain ray contractions of toric fans, which extends a key result by A.-M.~Castravet and J.~Tevelev \cite[Proposition~3.1]{castravet2015overline}. As a consequence, we are able to provide new instances of toric varieties whose blowup at a general point is not a Mori dream space and even does not have a polyhedral pseudoeffective cone (e.g. Lemma~\ref{lemma.Pdn.Losev.Manin_V2}). Based on this, we show the following.

\begin{theorem}      \label{thm.cone.pdn}
$\overline{P}_{d,n}$ is not a Mori dream space and does not have a polyhedral pseudoeffective cone 
for $n \geq d+9$.
\end{theorem}

\subsection{Strategy} Our strategy can be summarized as follows. We first study the case $X=\PP^d$. Dealing with $\PP^d[n]$ directly poses certain difficulties, for instance, all the centers in its iterated blowup construction have dimension larger than 0 and are not invariant under the torus of $(\PP^d)^n$. To overcome these issues, we consider weighted versions of $\PP^d[n]$, introduced by one of the authors in \cite{routis2014weighted}, as well as certain GIT quotients of them, introduced and studied by two of the authors in \cite{GR}. These quotients can be regarded as higher-dimensional generalizations of Hassett's moduli space of weighted rational stable curves \cite{Hassett-weighted}, as they satisfy several analogous properties, including the existence of weight reduction morphisms. In particular, they admit surjective morphisms from $\overline P_{d,n}$. For any $d$ and $n>d+2$, we are able to select a suitable set of weights and identify the associated quotient as the blowup of a toric variety $\overline P_{d,n}^{LM}$ at the unit of its torus $e$. Therefore, using results from \cite{okawa2016images} and \cite{baker2011good} we reduce the original question to showing that $Bl_e \overline P_{d,n}^{LM}$ is not a Mori dream space. We remark that a direct application of the key result \cite[Proposition~3.1]{castravet2015overline} to $\overline P_{d,n}^{LM}$ \emph{does not} seem to lead to a conclusive answer. Instead, our general method developed in Section~\ref{sec:technical} allows us to map this blowup, via a small modification followed by a surjection, to the blowup of $\overline P_{d-1,n-1}^{LM}$ at the unit of its torus. Applying this process repeatedly, we arrive at the blowup of some Losev-Manin space \cite{losev2000new-moduli-spaces} at the unit of its torus, from where we conclude the result by \cite{castravet2015overline}, \cite{gonzalez2016some}, \cite{hausen2018blowing}. Finally, Proposition~\ref{birFM} allows us to prove the more general case stated in Theorem~\ref{thm.FM}. \\

\subsection{Notation}
  Given a finite list $v_1,\linebreak[1]  v_2,\linebreak[1] \ldots,\linebreak[1] v_n$,\linebreak[1] for any $1 \leq j \leq n$, we denote the sublist obtained omitting $v_j$ by $v_1,\linebreak[1] v_2,\linebreak[1] \ldots,\linebreak[1] \widehat{v_{j}},\linebreak[1] \ldots ,\linebreak[1] v_n$. The cone generated by the vectors $\{ v_1,v_2, \ldots, v_n \}$ in some $\mathbb{R}^m$ is denoted as $\operatorname{Cone}\left(v_1,v_2, \ldots, v_n \right)$. \\
\indent Given an abelian group $N$ we use the notation $N_{\mathbb{R}}=N \otimes_\mathbb{Z} \mathbb{R}$ and given a homomorphism $\pi:N \rightarrow N'$ between free abelian groups, for simplicity we also denote by $\pi$ its scalar extension 
$\pi \otimes_\mathbb{Z} id_{\mathbb{R}}:N_{\mathbb{R}} \rightarrow N'_{\mathbb{R}}$.  \\ 
\indent We use the term variety to refer to a reduced and irreducible scheme of finite type over an algebraically closed field of characteristic 0. 
The group of Weil divisors of a variety $X$ is denoted as  $WDiv(X)$. 
A small $\mathbb{Q}$-factorial modification (SQM for short) of a normal, projective, $\mathbb{Q}$-factorial variety $X$ is a birational map $X \dashrightarrow Y$ to another normal, projective, $\QQ$-factorial variety $Y$, which is an isomorphism in codimension one.\\ 
\indent For any reductive group $G$ acting on a projective variety $X$ and any $G$-linearized line bundle $L$ on $X$, we denote the semistable locus in $X$ with respect to $L$ as $X^{ss}(L)$ and the stable locus as $X^s(L)$.

\subsection{Acknowledgements}
The authors would like to thank W. Fulton for insightful conversations. P. Gallardo thanks the University of California, Riverside for its support. J. Gonz\'alez was supported by a grant from the Simons Foundation (Award Number 710443). E. Routis was supported by the EPSRC grant `Classification, Computation, and Construction: New Methods in Geometry' (EP/N03189X/1).


\section{Background} \label{section.background}

\subsection{Weighted Fulton-MacPherson spaces}\label{FM}
Let $X$ be a smooth variety with $\dim X=d\geq1$ and
let $n\geq 2$. 
The Fulton-MacPherson space $X[n]$ 
provides a natural compactification 
of the configuration space of $n$ distinct labeled points in $X$ \cite{Fulton-MacPherson}. 
\begin{remark}
We remark that $X[n]$ is proper if and only if $X$ is proper. Still, in the literature $X[n]$ is commonly called a compactification even when $X$ is not proper. We will follow that convention throughout this article. 
\end{remark}
Our work makes use of a weighted generalization of the above space
introduced by one of the authors in \cite{routis2014weighted}. The \emph{domain of admissible weights} for 
such generalization is 
\begin{align*}
\mathcal D_{d,n}^{FM} :=
\left\{
(a_1, \ldots,  a_n) \in \mathbb{Q}^{n} \; : \;
0 < a_i\leq 1\,\, , \ i=1,\dots, n
\right\}.
\end{align*}
For any $\mathcal{A} \in \mathcal D_{d,n}^{FM}$, we set  
\begin{equation*}\label{diag}
\mathcal{K}_{\mathcal{A}}
:= \{ \Delta_I \subseteq X^n \; | \; I \subseteq \{1,\dots,n\} \,\text{and}\, \sum \limits _{i\in I} a_i >1\}
,  \,\text{where}\,\,\,
\Delta_I:=
\{ (x_1, \ldots, x_n) \in X^n | x_i=x_j, \; \text{for all}\,\, i,j \in I \}.
\end{equation*}
 
It follows that for each $\Delta_I\in \mathcal{K}_{\mathcal{A}}$ we have $|I|\geq 2$. We will refer to $\Delta_I$ as the \emph{diagonal} corresponding to $I$. 
\begin{definition}
The weighted compactification $\wmod$ of $X^n\setminus \bigcup\limits_ {\Delta_I \in \mathcal{K}_{\mathcal{A}}} \Delta_I$ 
is the closure of the image of the natural locally closed embedding
$$
X^n\setminus \bigcup\limits_ {\Delta_I \in \mathcal{K}_{\mathcal{A}}} \Delta_I
\longrightarrow
\prod_{
\Delta_I \in \mathcal{K}_{\mathcal{A}}
} Bl_{\Delta_I}X^n.
$$
 
\end{definition}

\begin{definition}[{\cite[Definition 2.7]{LiLi}}]  \label{dominant} Let $Z$ be a subvariety of a variety $Y$ and $\pi : Bl_ZY \rightarrow Y$ be the blow-up of $Y$ along $Z$. For any subvariety $V$ of $Y$, the \textit{dominant transform} of $V$ is the strict transform of $V$ if $V \not\subseteq Z$ or the scheme-theoretic inverse $\pi^{-1}(V )$ if $V \subseteq Z$. \end{definition}

\begin{definition}
 Let $\mathcal{G}=\{S_1,S_2,\dots ,S_n\}$ be a totally ordered set of subvarieties of a variety $Y$. For $k=0,1,\dots, n$, we define $(Y_k,\{S_1^{(k)},S_2^{(k)},\dots ,S_n^{(k)}\})$ inductively as follows.
\begin{itemize}[leftmargin=*]
\item Let $Y_0=Y, S_i^{(0)}=S_i$ for $1\leq i\leq n$.
\item Let $k\geq1$ and suppose that $(Y_{k-1}, \{S_1^{(k-1)},S_2^{(k-1)},\dots,S_n^{(k-1)}\} )$ has been defined. 
\begin{itemize}
\item We define $Y_k$ as the blowup of $Y_{k-1}$ at the subvariety $S_k^{(k-1)}$.
\item We define $S_i^{(k)}$ as the dominant transform of $S_i^{(k-1)}$ under the blowup $Y_k\rightarrow Y_{k-1}$ for $1\leq i\leq n$.
\end{itemize}
\end{itemize}
The \textit{iterated blowup} $Bl_{\mathcal{G}} Y$ of $Y$ at $\mathcal{G}$ is the variety $Y_n$. The \textit{iterated dominant transform} of $S_i$ in $Y_n$ is the variety $S_i^{(n)}$.
\end{definition}

\vspace{0.1 in} 

The main properties of $\wmod$ may be summarized as follows (see \cite[Theorems 2 and 3]{routis2014weighted}).
\begin{theorem}\label{wcomp} 
For any smooth variety $X$ with $\dim X=d\geq1$, any $n\geq 2$ and any $\mathcal A\in \mathcal D^{FM}_{d,n}$:
\begin{enumerate}
\item {$\wmod$ is a nonsingular variety. Further, if $\sum \limits _{i=1}^n a_i >1$, the boundary $X_\mathcal{A}[n] \setminus (X^n \setminus \bigcup \limits_ {\Delta_I \in \mathcal{K}_{\mathcal{A}}} \Delta_I) $ is the union of $|\mathcal{K}_{\mathcal{A}}|$ smooth irreducible divisors $D_I,$ where $I\subseteq \{1,\dots, n\} $ and $ \sum \limits _{i\in I} a_i >1\ $.}

\item Let $\mathcal{K}_{\mathcal{A}}$ be given a total order $(<)$ 
extending the partial order by ascending dimension, that is $\Delta_I<\Delta_J$ if $dim(\Delta_I) <dim(\Delta_J)$. Then $X_\mathcal{A}[n]$ is the iterated blowup of $X^n$ at $\mathcal{K}_{\mathcal{A}}$. Moreover, each divisor $D_I$ is the iterated dominant transform of $\Delta_I$ in $\wmod$.  
 \end{enumerate}
 \end{theorem} 

 In particular, when $\mathcal{A}= (1, \ldots, 1)$,  the variety $X_\mathcal{A}[n]$ coincides with the Fulton-MacPherson compactification $X[n]$ (see \cite[Section 4.2]{LiLi} for the construction of $X[n]$ in an order of increasing dimension of diagonals). 
 
\begin{remark}In the statement of Theorem~\ref{wcomp} (2), the resulting iterated blowup is independent of the (total) order it is performed as long as it preserves the partial ascending dimension order of the diagonals.

\end{remark}
 
\subsection{Wonderful compactifications of the moduli space of points
in projective space}      \label{subsection.wonderful}
Next, we describe certain GIT quotients of weighted Fulton-MacPherson spaces that play a central role in our work.

Let us fix $d$ and $n$ such that $d\geq 1$ and $n > d+2$ and let $\epsilon=\frac{1}{n-d}$ and  $\hat \epsilon=\frac{1}{(d+1)(n-d)}$. Let
\begin{align*}
w_1=\cdots =w_d=1-\hat \epsilon,
& &
w_{d+1}=1-(n-(d+1))\epsilon +d \hat \epsilon, 
& & 
 w_{d+2}=\cdots = w_n=\epsilon. 
\end{align*}
The domain of admissible weights for such quotients is
\begin{align*}\label{eq:WeightPdn} 
\mathcal D_{d,n}^P := \left\{ (a_1,  \ldots , a_n) \in \mathbb{Q}^{n} \; : 
w_i \leq a_i \leq 1\,\,,i=1,\dots, n
\right\}.
\end{align*}

Let $\mathcal A\in  \mathcal D_{d,n}^P$. By Theorem~\ref{wcomp}, $\bP^d_{\mathcal A}[n]$ is the iterated blowup of $(\bP^d)^n$ at $\mathcal{K}_{\mathcal{A}}=\{ \Delta_I \subseteq (\PP^d)^n \; | \; I \subseteq \{1,\dots,n\} \,\text{and}\, \sum \limits _{i\in I} a_i >1\}$ given any total order which respects the partial order by ascending dimension. Let $\pi_{\mathcal A}:\bP^d_{\mathcal A}[n] \to (\bP^d)^n$ be the resulting morphism, and 
let $E_{\mathcal{A}}$ be the sum of its exceptional divisors.

Since $SL_{d+1}$ has no non-trivial characters, each ample line bundle has a unique $SL_{d+1}$-lin\-e\-ar\-i\-za\-tion for the diagonal action of $SL_{d+1}$ on $(\bP^d)^n$ induced by the usual action of $SL_{d+1}$ on $\bP^d$. 		
Moreover, given an $SL_{d+1}$-linearized ample line bundle, the semistable loci of each of its powers are the same.	 	
The same holds for the stable loci. 
Consider the line bundle  $$L_{d,n}=\mathcal{O}_{(\bP^{d})^n}\bigl((d+1)(n-d)w_1,(d+1)(n-d)w_2,...,(d+1)(n-d)w_n\bigr).$$ We have the following.  

\begin{lemma}[{\cite[Lemma~4.2]{GR}}]  \label{GITpoints} 
There exists an isomorphism
\begin{equation*}\label{GITiso}(\bP^d)^n /\!\!/_{L_{d,n}} SL_{d+1}\iso \left( \bP^{n-d-2} \right)^{d}.\end{equation*}
Furthermore, there is no strictly semistable locus in $(\bP^d)^n$ with respect to $L_{d,n}$.
\end{lemma}

  The diagonal action of $SL_{d+1}$ on $(\bP^d)^n$ lifts to $\bP^d_{\mathcal A}[n]$ \cite[Lemma~3.15]{GR}. 
  By the results in \cite[Section 4.1]{GR}, we have that 
for all positive integers $e$ sufficiently large,
  $\tilde{L}_{\mathcal A,e}:=\pi_{\mathcal A}^*(L_{d,n}^{\otimes e}) \otimes \mathcal{O}(-E_{\mathcal{A}})$ is a very ample line bundle which admits a unique $SL_{d+1}$-linearization such that
\begin{equation}\label{stablelocus}
(\bP^d_{\mathcal A}[n] )^{ss}(\tilde{L}_{\mathcal A,e})=(\bP^d_{\mathcal A}[n] )^s(\tilde{L}_{\mathcal A,e}) = \pi_{\mathcal A}^{-1}[((\bP^d)^n)^s(L_{d,n})].
\end{equation}
Then, for all 
large enough $e$, 
the line bundles $\tilde{L}_{\mathcal A,e}$ yield the same GIT quotient. 
We may thus take the smallest such $e$, which we denote by $e_0$, and
set $\tilde L_{\mathcal A}:=\tilde{L}_{\mathcal A,e_0}$ to define the GIT quotient
 $$
\overline P^{\mathcal A}_{d,n}:= \bP^d_{\mathcal A}[n] / \! \!/_{\tilde L_{\mathcal A}} SL_{d+1}.
$$

\begin{remark}\label{geomquot}By (\ref{stablelocus}) it follows that the above GIT quotient is actually a geometric quotient. 
\end{remark}

\begin{remark} When $d=1$, the variety $\overline P^{\mathcal A}_{1,n}$ is identified with Hassett's moduli space of rational weighted stable curves $\overline{M}_{0, \mathcal{A}}$ \cite{Hassett-weighted}.
\end{remark}
\begin{definition}\label{def:linebundle}
When $\mathcal{A}=(1,1,\dots,1)$ we will simply write $\overline P_{d,n}$ in place of $\overline P^{\mathcal A}_{d,n}$ and $\tilde L_{d,n}$ in place $\tilde L_{\mathcal A}$.
\end{definition}
The properties of $\bP^d_{\mathcal A}[n]$ described in Theorem~\ref{wcomp} descend to the quotient $\overline P^{\mathcal A}_{d,n}$, see \cite[Section 4.1]{GR} for more details. Here we only record some of them, which will be useful in the sequel. The following definition is given in \cite[Definition 4]{GR}.
\begin{definition}\label{hi}
\begin{enumerate} \item Let $[b^k_{d+2} : \ldots: b^k_{n}], \, k=0,1,\dots, d-1,$ be a system of projective coordinates for each copy of ${\PP}^{n-d-2}$ in the product  $\left( \bP^{n-d-2} \right)^d$. Let $I\subsetneq \{d+1,\dots n\}$ such that $ |I| \geq 2$. We define subvarieties $H_I$ of $({\PP}^{n-d-2})^d$ as follows:
  \begin{equation*}
  H_I :=\begin{cases}  
\bigcap_{i,j \in I} V((b^0_i-b^0_j, \ldots ,b^{d-1}_i-b^{d-1}_j)),  &  \  \text{if } d+1 \notin I 
\\
\bigcap_{i\in I \setminus d+1} V((b^0_i,\ldots,b^{d-1}_i)),        &  \   \text{if } d+1 \in I. 
\end{cases}
  \end{equation*}
  \item Let $\mathcal{A}=\{a_1,a_2,\dots, a_n\}\in \mathcal D_{d,n}^P$. We define the set $$\mathcal{G}_{\mathcal{A}}:=\{ H_I  \subseteq ( \bP^{n-d-2})^d  \; | \;I\subsetneq \{d+1,\dots, n\}\,\, \text{and}\,\, \sum_{i \in I} a_i > 1 \}.$$
  \end{enumerate}
\end{definition}

We may now state the following.

\begin{proposition}[{\cite[Section 4.1]{GR}}]    \label{smoothproj}

\begin{enumerate}\item The subvarieties $H_I$ of $({\PP}^{n-d-2})^d$ parametrize those equivalence classes of configurations of $n$ points $p_1,\dots,p_n$ in $((\bP^d)^n)^s(L_{d,n})$ where all points $ p_{i}$ such that $i \in I $ coincide.  Furthermore, $H_I\iso  ({\PP}^{(n-|I|)-d-1})^d$.

\item Let $\mathcal{G}_{\mathcal{A}}$ be given any total order $(<)$ extending the partial order by ascending dimension,  that is $H_I<H_J$  if $dim(H_I)<dim(H_J)$. Then, $\overline P^{\mathcal A}_{d,n}$ is the iterated blowup of $\left( \bP^{n-d-2} \right)^d$ at $\mathcal{G}_{\mathcal{A}}$ in that order.
\item For any $\mathcal A\in \mathcal D^{P}_{d,n}$, the variety $\overline P^{\mathcal A}_{d,n}$ is smooth and projective.

\end{enumerate}
\end{proposition}

\subsubsection{Toric models $\overline P_{d,n}^{LM}$}\label{toric.pdn} Consider the set of weights
$\mathcal{A}=(a_1,\dots,a_n)\in\mathcal D_{d,n}^{P}$ where 
$$
a_1 = \cdots = a_{d+1} =1, \;\;
a_{d+2} = \cdots = a_n =\frac{1}{n-d-1}.
$$
For such $\mathcal{A}$, the variety $\overline P_{d,n}^{LM}:=\overline P^{\mathcal A}_{d,n}$ is a toric variety which can be viewed as a higher dimensional analog of the Losev-Manin space. In particular, $\overline P_{1,n}^{LM}$ coincides with the Losev-Manin space of dimension $n-3$ \cite{losev2000new-moduli-spaces}.

Let us describe its toric structure. Consider the toric variety $\PP^{n-d-2}$ whose fan has rays generated by the following points in the lattice $\ZZ^{n-d-2}$  $$e_{d+2}:=(1,0,\dots,0),\dots, e_{n-1}:=(0,0,\dots,1),e_n:=(-1,-1,\dots,-1).$$

For $d+2\leq k \leq n$ and $1\leq i \leq d$, we define $e_k^i$ to be the image of $e_k$ under the inclusion 
\begin{gather*} 
\ZZ^{n-d-2}\to (\ZZ^{n-d-2})^d\\
x\mapsto \underbrace{(0,\dots,x,\dots,0)}_{i-\text{th position}}.
\end{gather*}
The $e_k^i$ generate the rays of the fan of the toric variety $(\PP^{n-d-2})^d$.
\begin{proposition}[{\cite[Corollary~5.7]{GR}}]
$\overline P_{d,n}^{LM}$ is a toric variety whose fan has rays generated by $e^i_{d+2}, \ldots, e^i_{n}$, where $1\leq i \leq d$, and 
$
\sum_{i\in I}\left( e^1_{i} + \cdots +e^{d}_{i} \right)
$
where  $1 \leq |I| \leq n-d-2$  and $I \subseteq \{d+2, \ldots, n\}$.
\end{proposition}
 
\begin{example}
We discuss the case $d=2$ and $n=5$. 
By  Proposition~\ref{smoothproj}, for any $\mathcal A\in  \mathcal D_{2,5}^P$, the variety 
$\overline  P_{2,5}^{\mathcal{A}}$ is the blow up of $\bP^1 \times \bP^1$ along 
$\mathcal{G}_{\mathcal{A}}$. Using the notation of Definition 
\ref{hi}, we should consider the indices
$\{3,4\}$, $\{3,5\}$ and $\{ 4,5\}$ which yield the following  subvarieties of $\PP^1\times \PP^1$
\begin{align*}
H_{\{3, 4\}} = ([0:1],[0:1]),    
&&
H_{\{3, 5 \}} = ([1:0],[1:0]),    
&&
H_{\{4, 5 \}} = ([1:1],[1:1]). 
\end{align*}
Given the weights $\mathcal{A}_1:=(1, 1,1, 1, 1)$ and 
$\mathcal{A}_2:=(1,1, 1, \frac{1}{2}, \frac{1}{2})$, we obtain the sets
\begin{align*}
\mathcal{G}_{\mathcal{A}_1} 
&=
\left\lbrace
([0:1],[0:1]), \; ([1:0],[1:0]), \; ([1:1],[1:1])
\right\rbrace 
\\
\mathcal{G}_{\mathcal{A}_2}  &=
\left\lbrace
([0:1],[0:1]), \; ([1:0],[1:0])
\right\rbrace.
\end{align*}
Therefore, writing $e$ for $([1:1],[1:1])$, we have that 
\begin{align*}
\xymatrix{
\overline  P_{2,5} \cong \operatorname{Bl}_{e,2 pts} 
\left( \bP^1 \times \bP^1 \right)
&& 
\overline  P_{2,5}^{LM} \cong \operatorname{Bl}_{2 pts} 
\left( \bP^1 \times \bP^1 \right)
}.
\end{align*}
We arrive at the sequence of morphisms
\begin{align*}
\xymatrix{
\overline  P_{2,5}  \ar[r]^{\cong \qquad}
&
\operatorname{Bl}_{e}  \overline  P^{LM}_{2,5}  \ar[r]
&
\overline  P_{2,5}^{LM} \ar[rr]
&& 
\bP^1 \times \bP^1 
}
\end{align*}
which appropriately generalizes to all $d\geq 1$ and $n> d+2$, see Proposition~\ref{thm:SecblowUps} (the first arrow is not an isomorphism in general).
\end{example}

\subsection{Preliminaries on Mori Dream spaces}\label{sec:MDS}
Let $X$ be a normal variety with a finitely generated divisor class group and only constant invertible global regular functions. For a finitely generated semigroup $\Gamma$ of Weil divisors on $X$, consider the $\Gamma$-graded ring:
\begin{align*}
R\left( X, \Gamma \right)
:=
\bigoplus_{D \in \Gamma}
H^0\left( X, \mathcal{O}_X(D)  \right)   
\end{align*}
where $\mathcal{O}_X(D)$ is the divisorial sheaf associated to the Weil divisor $D$.  

If $\Gamma$ is a finitely generated  semigroup of Weil divisors on $X$ such that their classes span the vector space  $Cl(X)_{\mathbb{Q}}$, the ring $R(X, \Gamma)$ is called \emph{a Cox ring} of $X$. The definition of $R(X, \Gamma)$ depends on the choice of $\Gamma$, but whether it is finitely generated as an algebra over the base field does not depend on the choice of $\Gamma$. One frequently finds in the literature the \emph{Cox ring} defined as
$$
\mathcal{R}\left( X \right)
:=
\bigoplus_{ D \in Cl(X)} H^0\left(
X, \mathcal{O}_X(D)
\right).
$$
The finite generation of a Cox ring of $X$ is equivalent under either of the definitions presented above.  
Since our interest is in the failure of the finite generation of the Cox ring, we will move between these definitions at will. 

A normal, projective, $\mathbb{Q}$-factorial variety $X$ with a finitely generated class group is called a \emph{Mori dream space}, MDS for short, if it has a finitely generated Cox ring. 
A Cox ring of $X$ can be identified with a Cox ring for each of its small $\mathbb{Q}$-factorial modifications. Hence $X$ is an MDS if and only if any given small $\mathbb{Q}$-factorial modification of $X$ is an MDS. 
Moreover, the image of an MDS is an MDS, as shown by S.~Okawa \cite{okawa2016images}. 

\begin{theorem}[{\cite[Theorem~1.1]{okawa2016images}}]   \label{thm:PropertiesMDS}
Let $X$, $Y$ be normal, projective, $\mathbb{Q}$-factorial varieties, and $f : X \rightarrow Y$ be a surjective morphism. If $X$ is a Mori dream space, then so is $Y$.
\end{theorem}

H.~B\"aker showed that if the Cox ring of $X$ is finitely generated, the finite generation also holds for the sections  over open subsets of $X$ of the sheaves of algebras considered before. 

\begin{theorem}[{\cite[Theorem~1.2]{baker2011good}}] \label{thm:bak1B}
Let $X$ be a normal variety with finitely generated Cox ring $\mathcal{R}(X)$. Then, for any finitely generated subgroup $K \subseteq WDiv(X)$, and any open subset $U \subseteq X$, the algebra of sections $$\bigoplus_{D \in K}  \Gamma(U,  \mathcal{O}_X(D))$$ is finitely generated.
\end{theorem}

H.~B\"aker also showed that a good quotient of an MDS is again an MDS.  More precisely, the following result will be important in the sequel. 

\begin{theorem}[{\cite[Theorem~1.1]{baker2011good}}] \label{thm:bak1}
Let $X$ be a normal variety with finitely generated Cox ring $\mathcal{R}(X)$. Let $G$ be a reductive affine algebraic group acting on $X$, and let $U \subseteq X$ be an open invariant subset admitting a good quotient 
$$
\pi: U \to U / \! \! / G,
$$ 
that is, 
a $G$-invariant affine surjective morphism $\pi$ with  $\mathcal{O}_{U / \! \! / G} = \pi_*\left( \mathcal{O}_U \right)^G$ (this holds for a good quotient in the sense of 
\cite[Definition~1.5]{seshadri1972quotient}).  
If $U/ \!\! /G $ has only constant invertible global functions,  
    then
    $
    \mathcal{R}\left( 
    U/ \!\! /G 
    \right)$ is finitely generated.
\end{theorem}

Recall that $\overline P_{1,n}^{LM}$ is identified with the Losev-Manin compactification of dimension $(n-3)$, and let 
$e \in \overline{P}_{1,n}^{LM}$ be the unit of the torus. The following result first appeared in
\cite{castravet2015overline} for $n \geq 134$. The bound was then improved to $n \geq 13$ in \cite{gonzalez2016some}, and to $n \geq 10$ in \cite{hausen2018blowing}. 
It was used in the proof that $\overline{M}_{0,n}$ is not an MDS for $n\geq 10$ and is also important for our work. 
\begin{theorem} \label{LM}
$Bl_{e}\overline P_{1,n}^{LM}$ is not an MDS for $n \geq 10$.
\end{theorem}

\subsection{Preliminaries on polyhedrality of pseudoeffective cones}

Let $X$ be a normal, projective, $\mathbb{Q}$-factorial variety. 
Let $\operatorname{N}^1(X)$ be the group of numerical equivalence classes of Cartier divisors on $X$. 
It is well-known that the scalar extension $\operatorname{N}^1(X)_{\mathbb{R}}$ is a finite dimensional real vector space. 
The pseudoeffective cone of $X$ is the closure of the convex cone in $\operatorname{N}^1(X)_{\mathbb{R}}$ generated by the numerical classes of effective Cartier divisors. 
The pseudoeffective cone plays an important role in the study of the birational geometry of algebraic varieties (see for example \cite{LazarsfeldPositivity1} for further motivation). 

Recall that a convex cone $\sigma$ in $\mathbb{R}^n$ is called polyhedral if there exist finitely many vectors $v_1,v_2,\ldots,v_r \in \mathbb{R}^n$ such that $\sigma=Cone(v_1,v_2,\ldots,v_r)=\mathbb{R}_{\geq 0}v_1+\mathbb{R}_{\geq 0}v_2+\ldots+\mathbb{R}_{\geq 0}v_r$. 
The convex cone $\sigma$ is said to be rational polyhedral if additionally one can choose the $v_i$'s in $\mathbb{Q}^n \subseteq \mathbb{R}^n$.

If $Y$ is an SQM of $X$, the rational map $X\dashrightarrow Y$ induces a linear isomorphism between $\operatorname{N}^1(X)_{\mathbb{R}}$ and $\operatorname{N}^1(Y)_{\mathbb{R}}$ under which their pseudoeffective cones get identified. 
In particular, $X$ has a (rational) polyhedral pseudoeffective cone if and only if $Y$ does. As with the Mori dream space property, the (rational) polyhedrality of the effective cone also passes along a surjective morphism.

\begin{theorem}[{\cite[Lemma~2.2]{CastravetLafaceTevelevUgaglia2020}}]   \label{thm:PropertiesPseudoeffective}   
Let $X$, $Y$ be normal $\mathbb{Q}$-factorial projective varieties, and $f : X \rightarrow Y$ be a surjective morphism. If the pseudoeffective cone of $X$ is (rational) polyhedral, then the same holds for $Y$.
\end{theorem}

The pseudoeffective cone of a Mori dream space is rational polyhedral as it is generated by the classes of effective divisors associated to a set of generators of the Cox ring. 
The converse does not necessarily hold as there are non-Mori dream spaces with a rational polyhedral pseudoeffective cone.  
Adding to Theorem~\ref{LM} in this direction, A.-M.~Castravet, A.~Laface, J.~Tevelev and L.~Ugaglia showed the following. 

\begin{theorem}[{\cite[Proof of Theorem~1.3]{CastravetLafaceTevelevUgaglia2020}}] \label{LM_pseudoeffective}
The pseudoeffective cone of $Bl_{e}\overline P_{1,n}^{LM}$ is not  polyhedral for $n \geq 10$.
\end{theorem}


\section{Proof of Main Theorems}

In Section~\ref{sec:technical} we will establish some necessary results about factorizations of rational maps between certain blowups of toric varieties, which may also be of independent interest. We subsequently use these results to prove our main theorems in Sections \ref{proof.thm.1.2} and \ref{proof.thm.1.1}.

\subsection{A factorization of a rational map between blowups of toric varieties}
\label{sec:technical}

Let us consider a morphism $p:X \rightarrow X'$ of projective, $\mathbb{Q}$-factorial toric varieties arising from a surjective homomorphism of lattices $\pi: N \rightarrow N'$, and let $e$ and $e'$ be the units of the respective tori. 
In this subsection, we give a sufficient condition guaranteeing that $\operatorname{Bl}_{e'}X'$ is a Mori dream space or has a (rational) polyhedral pseudoeffective cone if $\operatorname{Bl}_{e}X$ has the corresponding property. 

\begin{lemma}      \label{lemma.small.modification}
Let $\widetilde{\pi_1}:\operatorname{Bl}_{(p,q)} (\mathbb{A}^m \times \mathbb{P}^n) \dashrightarrow \operatorname{Bl}_p \mathbb{A}^m$ be the rational map induced by the projection $\pi_1:~\mathbb{A}^m \times \mathbb{P}^n \rightarrow \mathbb{A}^m$, for some $p \in \mathbb{A}^m$ and $q \in \mathbb{P}^n$. 
Then, there exist a projective bundle $\psi:~P~\rightarrow~\operatorname{Bl}_p \mathbb{A}^m$ and a birational map $\theta: \operatorname{Bl}_{(p,q)} (\mathbb{A}^m \times \mathbb{P}^n) \dashrightarrow P$ such that $\widetilde{\pi_1} =\psi \circ \theta$ and such that $\theta$ is an isomorphism on the complement of some closed subsets of codimension at least two lying over $p~\in~\mathbb{A}^m$.
\end{lemma}

\begin{proof}
By applying a linear change of coordinates if necessary, we may assume that $p=0 \in \mathbb{A}^n$ and that $q$ is a torus invariant point of our choice. 
Let $x_1,x_2,\ldots,x_m$ and $y_1,y_2,\ldots,y_n$ denote the canonical bases of $\mathbb{R}^m$ and $\mathbb{R}^n$. 
Let us define $y_0=-y_1-y_2-\cdots-y_n$, $x=x_1+x_2+\cdots+x_m$ and $y=y_1+y_2+\cdots+y_n$. 
Given any $v \in \mathbb{R}^m$, let $\overline{v}=(v,0) \in \mathbb{R}^m \times \mathbb{R}^n$ and given any $w \in \mathbb{R}^n$, let $\overline{w}=(0,w) \in \mathbb{R}^m \times \mathbb{R}^n$.
We can consider $\mathbb{A}^m$ as the toric variety associated to the fan in $\mathbb{R}^m$ whose cones are the faces of $\operatorname{Cone}(x_1,x_2,\ldots,x_m)$ and 
we can consider $\mathbb{P}^n$ as the toric variety associated to the fan in $\mathbb{R}^n$ whose cones are the faces of $\operatorname{Cone}(y_0,y_1,\ldots,\widehat{y_j},\ldots,y_n)$, for $0 \leq j \leq n$.
We may assume that $q \in \mathbb{P}^n$ is the torus invariant point that corresponds to the maximal cone  $\operatorname{Cone}(y_1,y_2,\ldots,y_n)$. 

Then, $\operatorname{Bl}_p \mathbb{A}^m$ gets identified with the toric variety associated to the fan in $\mathbb{R}^m$ whose cones are all faces of $\operatorname{Cone}(x_1,x_2,\ldots,\widehat{x_i},\ldots,x_m,x)$, for each $1 \leq i \leq m$. 
Likewise, $\mathbb{A}^m \times \mathbb{P}^n$ is identified with the toric variety associated to the fan in $\mathbb{R}^m \times \mathbb{R}^n$ whose cones are all the faces of $\operatorname{Cone}(\overline{x_1},\overline{x_2}, \linebreak[1] \ldots \linebreak[1] , \linebreak[1] \overline{x_m}, \linebreak[1] \overline{y_0},  \linebreak[1] \overline{y_1}, \linebreak[1]  \ldots, \linebreak[1] \widehat{\overline{y_j}}, \linebreak[1] \ldots,\overline{y_n})$, for each $0 \leq j \leq n$. 
Then, $\operatorname{Bl}_{(p,q)} (\mathbb{A}^m \times \mathbb{P}^n)$ gets identified with the toric variety associated to the fan $\Delta_B$ in $\mathbb{R}^m \times \mathbb{R}^n$ whose cones are all faces of  
$\operatorname{Cone}(\overline{x_1},\overline{x_2},\ldots,\widehat{\overline{x_i}},\ldots,\overline{x_m},\overline{y_1},\ldots,\overline{y_n},\overline{x}+\overline{y})$ for $1 \leq i \leq m$,
together with all the faces of 
$\operatorname{Cone}(\overline{x_1}, \linebreak[1] \overline{x_2}, \linebreak[1] \ldots, \linebreak[1] \overline{x_m}, \linebreak[1] \overline{y_1}, \linebreak[1] \ldots, \linebreak[1] \widehat{\overline{y_j}}, \linebreak[1] \ldots, \linebreak[1] \overline{y_n}, \linebreak[1] \overline{x}+\overline{y})$ for $1 \leq j \leq n$,
and all the faces of
$\operatorname{Cone}(\overline{x_1}, \linebreak[1] \overline{x_2}, \linebreak[1] \ldots, \linebreak[1] \overline{x_m}, \linebreak[1] \overline{y_0}, \linebreak[1] \overline{y_1}, \linebreak[1] \ldots, \linebreak[1] \widehat{\overline{y_j}}, \linebreak[1] \ldots, \linebreak[1] \overline{y_n})$ for $1 \leq j \leq n$.

Let $D$ be the torus invariant divisor on $\operatorname{Bl}_p \mathbb{A}^m$ associated to the ray generated by $x$
and let $\mathcal{E}= \mathcal{O}_{\operatorname{Bl}_p \mathbb{A}^m} \oplus \mathcal{O}_{\operatorname{Bl}_p \mathbb{A}^m}(D)^{\oplus n}$.
Notice that $\mathcal{E}$ is a vector bundle over the toric variety $\operatorname{Bl}_p \mathbb{A}^m$ that decomposes as a direct sum of toric line bundles.  
Hence, by \cite[Proposition~7.3.3]{CLS}, the projectivization $\psi: P=\mathbb{P}(\mathcal{E}) \rightarrow \operatorname{Bl}_p \mathbb{A}^m$ is a toric variety and it can be obtained from the fan $\Delta_P$ in $\mathbb{R}^m \times \mathbb{R}^n$ whose cones are all faces of $\operatorname{Cone}(\overline{x_1},\overline{x_2},\ldots,\widehat{\overline{x_i}},\ldots,\overline{x_m},\overline{y_0},\overline{y_1},\ldots,\widehat{\overline{y_j}},\ldots,\overline{y_n},\overline{x}+\overline{y})$ for each $1 \leq i \leq m$ and $0 \leq j \leq n$.
The fan $\Delta_B$ of $\operatorname{Bl}_{(p,q)} (\mathbb{A}^m \times \mathbb{P}^n)$ and the fan $\Delta_P$ of $P$ on $\mathbb{R}^m \times \mathbb{R}^n$ have the same rays, and then we get an induced birational map $\theta: \operatorname{Bl}_{(p,q)} (\mathbb{A}^m \times \mathbb{P}^n) \dashrightarrow P$. 
By construction, we have that $\widetilde{\pi_1} =\psi \circ \theta$. 

Let us show that $\theta$ is an isomorphism on the complement of some closed subsets of 
$\operatorname{Bl}_{(p,q)} (\mathbb{A}^m \times \mathbb{P}^n)$ and $P$ of
codimension at least two lying over $p~\in~\mathbb{A}^m$. 
Since $\Delta_B$ and $\Delta_P$ have the same set of rays,  
by \cite[Lemma~3.3.21]{CLS} it is enough to show that any cone that is in one of the fans $\Delta_B$ and $\Delta_P$ but not in the other one has a projection to $\mathbb{R}^m$ which intersects the interior of $\operatorname{Cone}(x_1,x_2,\ldots,x_m)$.
For this, it is enough to show that each cone in one of the fans $\Delta_B$ or $\Delta_P$ which projects to $\mathbb{R}^m$ inside the boundary of $\operatorname{Cone}(x_1,x_2,\ldots,x_m)$ is also in the other fan. 
Notice that all $\overline{x_i}$ and $\overline{x}$ project respectively to $x_i$ and $x$ in $\mathbb{R}^m$, and that all $\overline{y_j}$ and $\overline{y}$ project to zero in $\mathbb{R}^m$.  
Hence, it is enough to show that each cone in one of the fans $\Delta_B$ and $\Delta_P$ whose set of generators does not include $\overline{x}+\overline{y}$ and also does not simultaneously include all $\overline{x_i}$ for $1 \leq i \leq m$ is also in the other one of $\Delta_B$ and $\Delta_P$. 
But this last statement clearly holds, as we can verify from the explicit description of the cones in  $\Delta_B$ and $\Delta_P$.
\end{proof}

The following proposition is a generalization of \cite[Proposition~3.1]{castravet2015overline}, which corresponds to the case $n=1$ of the present statement.
In the proof we follow a similar strategy and construct the small modification needed for the extension to the higher dimensional case.

\begin{proposition}   \label{proposition.projection} 

Let $\pi: N \rightarrow N'$ be a surjective homomorphism of lattices and let $K$ be its kernel.  
Let $n$ be the dimension of $K$ and let $a_1,a_2,\ldots,a_{n}$ be a basis of $K$. 
Let $a_0=-a_1-a_2-\cdots-a_{n}$ and $\Gamma_K=\{a_0,a_1,\ldots,a_{n}\}$. 
Let $\Gamma = \{y_1,\ldots,y_{r}\}$ be a collection of distinct primitive elements of $N$, such that none of them is in $K$.
For each $1 \leq i \leq r$ and $0 \leq j \leq n$, let $\overline{y_i}$ and $\overline{a_j}$ denote the rays in $N_{\mathbb{R}}$ spanned by $y_i$ and $a_j$, respectively.  
Suppose that there exists a complete simplicial fan $\Delta'$ in $N'_{\mathbb{R}}$ whose rays are precisely $\{ \pi({\overline{y_1}}), \pi({\overline{y_2}}), \ldots, \pi({\overline{y_r}})  \}$, where there can be repetitions. Suppose also that the associated toric variety $X'=X(\Delta')$ is projective. Then, 

\textbf{(a)} There exists a complete, simplicial fan $\Delta$ in $N_{\mathbb{R}}$ whose distinct rays are precisely $\overline{y_i}$ and $\overline{a_j}$, for $1 \leq i \leq r$ and $0 \leq j \leq n$, such that the associated toric variety $X=X(\Delta)$ is projective and such that $\pi$ maps each cone of $\Delta$ onto a cone of $\Delta'$. In particular, the induced rational map $p: X \dashrightarrow X'$ is a surjective morphism. 

\textbf{(b)} Fix $p: X=X(\Delta) \rightarrow X'=X(\Delta')$ as in (a). 
Let $e \in X$ and $e' \in X'$ be the unit elements of their tori. 
Let $\widetilde{p}:\operatorname{Bl}_{e}X \dashrightarrow \operatorname{Bl}_{e'}X'$ be the induced rational map. 
There exists a normal, projective, $\mathbb{Q}$-factorial variety $Z$ and a rational map $\lambda: \operatorname{Bl}_{e}X \dashrightarrow Z$ that is an isomorphism in codimension one, such that the rational map $\mu= \widetilde{p} \circ \lambda^{-1}: Z \dashrightarrow \operatorname{Bl}_{e'}X'$ is a  surjective  morphism. 

\textbf{(c)} Fix $p: X=X(\Delta) \rightarrow X'=X(\Delta')$ as in (a). 
Let $e \in X$ and $e' \in X'$ be the unit elements of their tori. 
If $\operatorname{Bl}_{e}X$ is a Mori dream space then $\operatorname{Bl}_{e'}X'$ is also a Mori dream space.
If the pseudoeffective cone of $\operatorname{Bl}_{e}X$ is (rational) polyhedral, then the pseudoeffective cone of $\operatorname{Bl}_{e'}X'$ is also (rational) polyhedral.
\end{proposition}


\begin{proof}


\textbf{(a)} We proceed by induction on $q=r-|\{ \pi({\overline{y_1}}), \pi({\overline{y_2}}), \ldots, \pi({\overline{y_r}})  \}|$. 
Let $d$ be the dimension of $X'$. 
If $q=0$, we let $\Delta$ be the set of cones in $N_{\mathbb{R}}$ closed under taking faces and whose maximal cones are the cones generated by any collection of rays of the form $\overline{y_{i_1}}, \overline{y_{i_2}}, \ldots, \overline{y_{i_d}},\linebreak[1] \overline{a_{0}},\linebreak[1] \overline{a_{1}},\linebreak[1] \ldots,\linebreak[1] \widehat{ \overline{a_{j}} },  \ldots,\linebreak[1] \overline{a_n}$, where $\pi(\overline{y_{i_1}}), \linebreak[1]  \pi(\overline{y_{i_2}}),  \linebreak[1] \ldots,  \linebreak[1] \pi(\overline{y_{i_d}})$
span a maximal cone in $\Delta'$ and $0 \leq j \leq n$.  
Since $\Delta'$ is a complete, simplicial fan in $N'_{\mathbb{R}}$, it follows that $\Delta$ is a complete, simplicial fan in $N_{\mathbb{R}}$. 
Notice that each cone in $\Delta'$ maps onto a cone in $\Delta$, and in particular the rational map $p: X=X(\Delta) \dashrightarrow X'=X(\Delta')$ is a surjective morphism. 
Notice also that $p$ is proper since $\pi^{-1}(|\Delta'|)=|\Delta|$.

We now show that the complete toric variety $X$ is projective.  
Let $T$ denote the torus of $X$.  
Let $D_0$ be the $T$-invariant divisor on $X$ corresponding to the ray $\overline{a_0}$. 
Fix an ample Cartier divisor $A$ on $X'$ and fix $m \in \mathbb{Z}^+$ sufficiently large that for any $T$-invariant curve in $X$ such that $p(C)$ is a curve, we have that $D_0 \cdot C + m p_*(C) \cdot A > 0$.
It is enough to show that $D=D_0+mp^*{A}$ is ample. 
By the Toric Kleiman Criterion \cite[Theorem~6.3.13]{CLS} it is enough to show that $D \cdot C >0$, for each $T$-invariant curve $C$ in $X$. 
Such a curve $C$ corresponds to a codimension-one cone $\tau$ in $\Delta$. 
For any such $\tau$ there exists a cone in $\Delta'$ spanned by distinct rays of the form $\pi(\overline{y_{i_1}}), \pi(\overline{y_{i_2}}), \ldots, \pi(\overline{y_{i_d}})$, such that the spanning rays of $\tau$ have one of the four possible forms: 
(1) $\overline{y_{i_1}}, \overline{y_{i_2}}, \ldots, \overline{y_{i_{d-1}}},\overline{a_{1}}, \overline{a_{2}}, \ldots, \overline{a_{n}}$;
(2) $\overline{y_{i_1}}, \overline{y_{i_2}}, \ldots, \overline{y_{i_{d}}},\overline{a_{1}}, \overline{a_{2}}, \ldots, \widehat{ \overline{a_{j}} },\ldots,\overline{a_{n}}$, for some $j \neq 0$;
(3) $\overline{y_{i_1}}, \linebreak[1]  \overline{y_{i_2}}, \linebreak[1]  \ldots, \linebreak[1]  \overline{y_{i_{d-1}}}, \linebreak[1] \overline{a_{0}}, \linebreak[1]  \overline{a_{1}}, \linebreak[1]  \ldots, \linebreak[1]  \widehat{ \overline{a_{j}} }, \linebreak[1] \ldots, \linebreak[1]  \overline{a_{n}}$, for some $j \neq 0$;
(4) $\overline{y_{i_1}},  \linebreak[1] \overline{y_{i_2}},  \linebreak[1] \ldots,   \linebreak[1] \overline{y_{i_{d}}}, \linebreak[1]  \overline{a_{0}}, \linebreak[1] \overline{a_{1}}, \linebreak[1] \ldots, \linebreak[1] \widehat{ \overline{a_{j}} }, \linebreak[1] \ldots,  \linebreak[1] \widehat{ \overline{a_{l}} },\ldots,\overline{a_{n}}$, for some $0 < j < l \leq n$.

In case (1), $D_0$ does not intersect $C$ and $p(C)$ is a curve in $X'$, hence $D \cdot C = m p^*A \cdot C > 0$. 
In case (2), $p(C)$ is a point in $X'$ and if we denote by $\sigma$ the cone in $\Delta$ generated by $\tau$ together with $a_0$, then $D \cdot C = D_0 \cdot C = \frac{\operatorname{mult(\tau)}}{\operatorname{mult(\sigma)}} >0$, by \cite[Lemma~6.4.2]{CLS}.
In case (3), $p(C)$ is a curve in $X'$ and hence using the projection formula we get $D \cdot C = D_0 \cdot C + m p^*A \cdot C = D_0 \cdot C + m A \cdot p_*C > 0$, by the choice of $m$. 
In case (4), first notice that $p(C)$ is a point in $X'$ and that the only $T$-invariant divisors on $X$ intersecting $C$ are those associated to the rays $\overline{y_{i_1}}, \overline{y_{i_2}}, \ldots, \overline{y_{i_{d}}},\overline{a_{0}}, \overline{a_{1}}, \ldots, \overline{a_{n}}$.
We can choose $u \in Hom(N, \mathbb{Z})$ such that $u(v)=0$ for $v=y_{i_1}, y_{i_2}, \ldots, y_{i_d},a_1,a_2,\ldots,\widehat{a_j},\ldots,a_n$ and $u(a_0)>0$. 
Since $a_0+a_1+\cdots+a_n=0$, then $u(a_0)=-u(a_j)$. 
Let $D_j$ the $T$-divisor in $X$ associated to the ray $\overline{a_{j}}$ and let $\sigma$ be the cone in $\Delta$ spanned by $\overline{a_{j}}$ and $\tau$.  
We have that $D_j \cdot C = \frac{\operatorname{mult(\tau)}}{\operatorname{mult(\sigma)}} >0$.
Notice that $\operatorname{div}(\chi^u)=u(a_0)D_0-u(a_0)D_j+F$, where $F$ is a sum of $T$-divisors that do not intersect $C$.
In particular, $0= \operatorname{div}(\chi^u) \cdot C = u(a_0)D_0 \cdot C - u(a_0)D_j \cdot C$, and hence $D_0 \cdot C = D_j \cdot C > 0$. 
Putting this together, in case (4) we have $D \cdot C = D_0 \cdot C > 0$. Therefore, $A$ is ample and then $X$ is projective.



Notice that the fan $\Delta$ that we constructed for $q=0$ satisfies that the preimage under $\pi$ of any ray in $\Delta'$ is a union of cones in $\Delta$. Let us include this additional property as part of our induction argument.   
For the inductive step, assume that $q \geq 1$ and by reordering the rays of $\Delta$, if necessary, we may assume that $\pi({\overline{y_r}}) \in \{ \pi({\overline{y_1}}), \pi({\overline{y_2}}), \ldots, \pi({\overline{y_{r-1}}})  \}$.  
By induction hypothesis, there exists a projective, simplicial fan $\Delta_0$ whose distinct rays are precisely $\overline{y_i}$ and $\overline{a_j}$, for $1 \leq i \leq r-1$ and $0 \leq j \leq n$, satisfying the conditions in (a) with respect to $\Delta'$. 
We may additionally assume that for each $1 \leq i \leq r-1$ the set $\pi^{-1}(\pi({\overline{y_i}}))$ is a union of cones in $\Delta_0$. 
We define $\Delta$ to be the star subdivision of $\Delta_0$ along the vector $y_r$. 
Then, $\Delta$ is a complete, simplicial fan. We get an induced surjective morphism $X=X(\Delta) \rightarrow X(\Delta_0)$ which is projective by \cite[Proposition~11.1.6]{CLS}. Then, $X$ is projective and the induced rational map $p: X \dashrightarrow X'$ is a surjective morphism.
Let $\sigma$ be a cone in $\Delta$. 
If $\sigma$ is also in $\Delta_0$, then it is mapped by $\pi$ onto a cone in $\Delta'$. 
If $\sigma$ is not in $\Delta_0$, then by the definition of the star subdivision, $\overline{y_r}$ is one of the rays of $\sigma$. 
Moreover, there exist cones $\tau$ and $\sigma_0$ in $\Delta_0$, such that $y_r \notin \tau$, $y_r \in \sigma_0$, $\tau \subseteq \sigma_0$
and such that $\sigma$ is the cone spanned by $\tau$ and $\overline{y_r}$. 
The intersection of $\sigma_0$ with $\pi^{-1}(\pi(\overline{y_r}))$ is not empty, and hence it is a nonempty union of cones whose rays map under $\pi$ to $\pi(\overline{y_r})$. Let $\overline{y_j}$ be a ray of $\sigma_0$ such that $\pi(\overline{y_j})=\pi(\overline{y_r})$.
Since $\sigma_0$ is simplicial, $\tau$ and $\overline{y_j}$ span a face $\tau'$ of $\sigma_0$. 
By construction $\pi(\sigma)=\pi(\tau')$, and $\pi(\tau')$ is equal to a cone in $\Delta'$ since $\tau' \in \Delta_0$. 
Then, each cone in $\Delta$ is mapped by $\pi$ onto a cone in $\Delta'$. 
Notice that for each $1 \leq j \leq r-1$, if $\pi^{-1}(\pi(\overline{y_j})) \neq \pi^{-1}(\pi(\overline{y_r}))$, then the cones in $\pi^{-1}(\pi(\overline{y_j}))$ remained the same when passing from $\Delta_0$ to $\Delta$. Moreover, since $\pi^{-1}(\pi(\overline{y_r}))$ is a union of cones in $\Delta_0$, then after the star subdivision along $y_r$,  $\pi^{-1}(\pi(\overline{y_r}))$ is a union of cones in $\Delta$.
Then (a) follows by induction.


\textbf{(b)} 
Let $b:\operatorname{Bl}_{e}X \rightarrow X$ and $b':\operatorname{Bl}_{e'}X' \rightarrow X'$ be the blowup morphisms and let us denote their respective exceptional divisors by $E$ and $E'$.
Let $T'$ denote the torus of $X'$. 
Let us pullback the open cover $\{X' \smallsetminus \{e'\}\} \cup \{T'\}$ of $X'$ by the morphism $p \circ b: \operatorname{Bl}_{e}X \rightarrow X'$. 
We get an open cover $\{U\} \cup \{V\}$ of $\operatorname{Bl}_{e}X$, with $U=b^{-1}(p^{-1}(X' \smallsetminus \{ e' \}))$ and $V=b^{-1}(p^{-1}(T'))$. 
Notice that the restriction of $b$ gives an isomorphism between $U$ and $p^{-1}(X' \smallsetminus \{ e' \})$. 
The intersection $U \cap V$ can be identified via $b$ with $p^{-1}(T' \smallsetminus \{e'\})$, since $p(e)=e'$ and $U \cap V$ is the preimage under $p \circ b$ of the intersection $(X' \smallsetminus \{e'\}) \cap T'= T' \smallsetminus \{e'\}$.
We will now construct varieties, morphisms and rational maps as in the following commutative diagram.


\begin{equation*}
       \begin{tikzcd}[row sep=12mm, column sep = 8mm]
   U \cap \widetilde{V} \arrow[rr,leftrightarrow,"\sim"] \arrow[dr,hook] \arrow[dd,swap,hook] &&
   U \cap V  \arrow[rr, two heads] \arrow[dd,hook] \arrow[dr,hook] &&
   T'\smallsetminus \{e'\}  \arrow[dd,hook] \arrow[dr,hook]   
   \\
   & {\ \ \ \widetilde{V} \ \ \ }  
   &&     V  \arrow[from=ll, crossing over,dashed,leftarrow,"\theta" near start]   &&
   \operatorname{Bl}_{e'}T'   \arrow[from=ll, crossing over,dashed,"\widetilde{\pi_1} \quad"]  \arrow[dd,hook]  
   \\
   U \arrow[rr,leftrightarrow,near start,"="] \arrow[dr,hook] 
   && U \arrow[rr, two heads] \arrow[dr,hook]  && X'\smallsetminus \{e'\} \arrow[dr,hook] 
   \\
   & {\ \ \ Z \ \ \ } 
   \arrow[rr,dashed,leftarrow,"\lambda"] \arrow[from=uu, crossing over,hook] 
   &&  \operatorname{Bl}_{e}X = U \cup V \arrow[rr,dashed] \arrow[from=uu, crossing over,hook] &&   \operatorname{Bl}_{e'}X'  
   \arrow[from=2-2,to=2-6, bend right=18, crossing over,"\psi" near end, two heads]
   \arrow[from=llll, bend right =18, crossing over,"\mu" near end, two heads]
  \end{tikzcd}
\end{equation*}


By our assumption on the map of fans $\Delta \rightarrow \Delta'$, we see that $p|:p^{-1}(T') \rightarrow T'$ is a trivial fibration which can be naturally identified with $T' \times \mathbb{P}^n \rightarrow T'$. 
Let $q \in \mathbb{P}^n$ be the point such that $(e',q) \in T' \times \mathbb{P}^n$ corresponds to $e \in p^{-1}(T') \subseteq X$. 
We get an induced identification between $V=b^{-1}(p^{-1}(T'))=\operatorname{Bl}_{e} p^{-1}(T')$ and $\operatorname{Bl}_{(e',q)} T' \times \mathbb{P}^n$. 
Using this identification, the projection $\pi_1: T' \times \mathbb{P}^n \rightarrow T'$ induces a rational map $\widetilde{\pi_1}: V=\operatorname{Bl}_{e} p^{-1}(T') \dashrightarrow \operatorname{Bl}_{e'} T'$. 
Now we consider $T'$ as an open subset of $\mathbb{A}^m$, for $m=\operatorname{dim}X'$, and we restrict the construction in Lemma~\ref{lemma.small.modification} to this setting.
We deduce that there exist a projective bundle $\psi:~\widetilde{V}~\rightarrow~\operatorname{Bl}_{e'} T'$ and a birational map $\theta: V \dashrightarrow \widetilde{V}$ such that $\widetilde{\pi_1} =\psi \circ \theta$ and such that $\theta$ is an isomorphism on the complement of some closed subsets of codimension at least two lying over $e'~\in~T'$.
Since $U \cap V$ lies over $T' \smallsetminus \{e' \}$, then $\theta$ maps $U \cap V$ isomorphically onto $\theta(U \cap V)$.
Let $Z$ be the scheme obtained by gluing $U$ and $\widetilde{V}$ along the open subsets $U \cap V$ of $U$ and $\theta(U \cap V)$ of $\widetilde{V}$ which are identified via the isomorphism $\theta|_{U \cap V}$. 
By construction the induced rational map $\lambda: \operatorname{Bl}_{e}X \dashrightarrow Z$ is an isomorphism in codimension one. 
Let $\mu_U:U\rightarrow \operatorname{Bl}_{e'} X'$ obtained by composing the morphism $p \circ b:U \rightarrow X' \smallsetminus \{ e' \}$ with the inclusion $X' \smallsetminus \{ e' \} \subseteq \operatorname{Bl}_{e'}X'$.
Let $\mu_{\widetilde{V}}: \widetilde{V} \rightarrow \operatorname{Bl}_{e'} X'$ obtained by composing the morphism $\psi: \widetilde{V} \rightarrow \operatorname{Bl}_{e'}T'$ with the inclusion $\operatorname{Bl}_{e'}T' \subseteq \operatorname{Bl}_{e'}X'$.
Since $\widetilde{\pi_1} =\psi \circ \theta$, it follows that over $U \cap V$ both $(\mu_U)|_{U \cap V}$ and $(\mu_{\widetilde{V}} \circ \theta)|_{U \cap V}$ agree with $(\widetilde{\pi_1})|_{U \cap V}$ composed with the inclusion $\operatorname{Bl}_{e'} T' \subseteq \operatorname{Bl}_{e'} X'$.
Therefore, we can glue $\mu_U$ and $\mu_{\widetilde{V}}$ to get a morphism $\mu: Z \rightarrow \operatorname{Bl}_{e'} X'$ such that $\mu|_U=\mu_U$ and   $\mu|_{\widetilde{V}}=\mu_{\widetilde{V}}$.  
Since $\psi= \widetilde{\pi_1} \circ \theta^{-1}$ as rational maps, then $\mu= \widetilde{p} \circ \lambda^{-1}$.

We will now show that $Z$ is a normal, projective,  $\mathbb{Q}$-factorial  variety. 
Notice that $\mu^{-1}(\operatorname{Bl}_{e'} T')=\widetilde{V}$ and that $\mu|_{\widetilde{V}}=\psi:\widetilde{V}\rightarrow \operatorname{Bl}_{e'} T'$ is proper since it is a projective bundle. 
Notice also that $\mu^{-1}(X' \smallsetminus \{ e' \})=U=b^{-1}(p^{-1}(X' \smallsetminus \{ e' \}))$ and that 
$\mu_{U}=(p \circ b)|_{U}: U \rightarrow X' \smallsetminus \{ e' \}$ is proper since $p \circ b$ is proper, because both $p$ and $b$ are proper. 
Then, $\mu: Z \rightarrow \operatorname{Bl}_{e'} X'$ is a proper surjective morphism. 
In particular, $Z$ is proper since $\operatorname{Bl}_{e'} X'$ is proper.  
Since $Z$ is proper (in particular separated) and has an open cover by the normal,  $\mathbb{Q}$-factorial varieties $U$ and $\widetilde{V}$, we conclude that $Z$ is a normal, complete,  $\mathbb{Q}$-factorial variety. 
Since $\mu: Z \rightarrow \operatorname{Bl}_{e'} X'$ is a proper morphism and $\operatorname{Bl}_{e'} X'$ is projective, to show that $Z$ is also projective it is enough to show that $Z$ admits a $\mu$-ample Cartier divisor.
Let $A$ be the strict transform on $Z$ of an irreducible very ample Cartier divisor $A_0$ on the projective $\mathbb{Q}$-factorial variety $\operatorname{Bl}_e X$. 
Since $U= \mu^{-1}(T' \smallsetminus \{ e' \})$ and $A|_U=A_{0}|_U$, then $A|_{U}$ is $(\mu|_{U})$-ample. 
Notice that $A$ is mapped surjectively onto $\operatorname{Bl}_{e'} X'$ by $\mu$ since its image contains $T' \smallsetminus \{ e' \}$ and is closed by the properness of $\mu$. 
Since $\widetilde{V}= \mu^{-1}(\operatorname{Bl}_{e'} T')$ and $\mu|_{\operatorname{Bl}_{e'} T'}: \widetilde{V} \rightarrow \operatorname{Bl}_{e'} T'$ is a projective bundle, the irreducible Cartier divisor $A|_{\operatorname{Bl}_{e'} T'}$ is $(\mu|_{\operatorname{Bl}_{e'} T'})$-ample since it maps surjectively onto the base $\operatorname{Bl}_{e'} T'$. %
Therefore, $A$ is $\mu$-ample, hence $\mu$ is a projective morphism and $Z$ is projective. 

\textbf{(c)} Fix a variety $Z$ as in part $\operatorname{(b)}$. 
If $\operatorname{Bl}_{e}X$ is a Mori dream space or it has a (rational) polyhedral pseudoeffective cone, then $Z$ has the same respective property because they are normal, projective, $\mathbb{Q}$-factorial varieties isomorphic in codimension one. 
Now, if $Z$ has either of these properties, then the same holds for $\operatorname{Bl}_{e'}X'$ because $\mu:Z \rightarrow \operatorname{Bl}_{e'}X'$ is a surjective morphism of normal, projective, $\mathbb{Q}$-factorial varieties. 
Indeed, in the case where it has the Mori dream space property this holds by 
Theorem~\ref{thm:PropertiesMDS} and in the case where its pseudoeffective cone is (rational) polyhedral  this holds by Theorem~\ref{thm:PropertiesPseudoeffective}. 
\qedhere


\end{proof}


The following corollary provides the sufficient condition stated as the goal of this subsection.

\begin{corollary}     \label{corollary.projection}
Let $f:X(\Delta_1) \dashrightarrow X(\Delta_2)$ be a rational map between projective, $\mathbb{Q}$-factorial toric varieties induced by a surjective homomorphism of lattices $\pi: N_1 \rightarrow N_2$, and let $e_1$ and $e_2$ denote the units of the respective tori. Assume that: 
\begin{enumerate} 
    \item Each ray in $\Delta_1$ is mapped by $\pi$ to either zero or to a ray of $\Delta_2$, with possible repetitions. 
    \item The number of rays of $\Delta_1$ mapped to zero by $\pi$ is $\operatorname{dim}(X(\Delta_1))-\operatorname{dim}(X(\Delta_2))+1$, and the primitive generators of those rays add to zero and generate the kernel of $\pi$ as a group. 
\end{enumerate}
If $\operatorname{Bl}_{e_1}X(\Delta_1)$ is a Mori dream space,  then $\operatorname{Bl}_{e_2}X(\Delta_2)$ is also a Mori dream space. 
If the pseudoeffective cone of $\operatorname{Bl}_{e_1}X(\Delta_1)$ is (rational) polyhedral,  then the pseudoeffective cone of $\operatorname{Bl}_{e_2}X(\Delta_2)$ is also (rational) polyhedral. 
\end{corollary}
\begin{proof}
We apply Proposition~\ref{proposition.projection} to the homomorphism $\pi$, the fan $\Delta'=\Delta_2$, the set $\Gamma_K$ equal to the primitive generators of the rays of $\Delta_1$ in the kernel of $\pi$ and the set $\Gamma$ equal to the primitive generators of the remaining rays of $\Delta_1$. 
Then, $X'=X(\Delta')=X(\Delta_2)$ and notice that the rays of $\Delta_1$ that do not map to zero under $\pi$ map precisely to the rays of $\Delta'$. 
By Proposition~\ref{proposition.projection}, there exist (i) a toric variety $X=X(\Delta)$ such that $\Delta$ has the same rays as $\Delta_1$, and (ii) a normal, projective, $\mathbb{Q}$-factorial variety $Z$ that is a small $\mathbb{Q}$-factorial modification of $\operatorname{Bl}_e X$ and which admits a surjective morphism $Z \rightarrow \operatorname{Bl}_{e'} X'$, where $e'$ and $e$ denote the respective unit elements of the tori. 
Clearly, $\operatorname{Bl}_{e_1} X(\Delta_1)$, $\operatorname{Bl}_e X$ and $Z$ are normal projective $\mathbb{Q}$-factorial varieties isomorphic in codimension one. 
If $\operatorname{Bl}_{e_1} X(\Delta_1)$ is a Mori dream space or has a (rational) polyhedral pseudoeffective cone then the same holds for its SQM $Z$. %
Since $Z$ maps surjectively onto $\operatorname{Bl}_{e'} X'=\operatorname{Bl}_{e_2} X(\Delta_2)$, 
if $Z$ is a Mori dream space or has a (rational) polyhedral pseudoeffective cone then the same respective property holds for $\operatorname{Bl}_{e_2} X(\Delta_2)$. 
In the case of the Mori dream space property this holds by  Theorem~\ref{thm:PropertiesMDS} and in the case of the (rational) polyhedrality of the pseudoeffective cone this holds by  Theorem~\ref{thm:PropertiesPseudoeffective}.  
\qedhere
\end{proof}

\subsection{Proof of Theorem~\ref{thm.cone.pdn}}        \label{proof.thm.1.2}

Throughout this section, we fix a positive integer $d$. Recall from Section~\ref{subsection.wonderful} that $\overline{P}_{d,n}^{LM}$ is a toric variety. Let $e$ denote the unit element of its torus. We first show the following.


\begin{lemma}    \label{lemma.Pdn.Losev.Manin_V2} 
$\operatorname{Bl}_{e}\overline{P}_{d,n}^{LM}$ is not a Mori dream space for $n \geq d+9$  and moreover its pseudoeffective cone is not polyhedral.
\end{lemma}
\begin{proof} We proceed by induction on $d$. When $d=1$,  $\overline{P}_{d,n}^{LM}$ is the Losev-Manin compactification of dimension $n-3$. 
For $n\geq 10$, the variety $\operatorname{Bl}_{e}\overline{P}_{d,n}^{LM}$  is not an MDS by Theorem~\ref{LM} and its pseudoeffective cone is not polyhedral by Theorem~\ref{LM_pseudoeffective}. 

Now let $d\geq 2$ and suppose that $\operatorname{Bl}_{e}\overline{P}_{d-1,n-1}^{LM}$ is not a Mori dream space and its pseudoeffective cone is not polyhedral
for $n-1 \geq (d-1)+9$. 
Let us fix $n \geq d +9$.
Let $N_{d,n}$ be the lattice $(\ZZ^{n-d-2})^d$ with the points
$\{  e^i_{d+2}, \ldots, e^i_{n}\}$
defined in Section \ref{toric.pdn},
where $1 \leq i \leq d$. Let us define a homomorphism
 $$\pi_{d,n}:N_{d,n}\to N_{d-1,n-1}$$ 
\begin{equation*}
 e^i_k\mapsto
\begin{cases}
0, & \text{if}\,\,\, i=d\\
e^i_{k-1}, & \text{otherwise,}
\end{cases}
\end{equation*}
where  $d+2\leq k\leq n$.
\noindent Recall that the rays of the fan of the toric variety $\overline{P}_{d,n}^{LM}$ are spanned by the lattice points $e^i_{d+2}, \ldots, e^i_{n}$ and 
$
\sum_{k \in I}\left( e^1_{k} + \ldots +e^{d}_{k} \right)
$
where  $1 \leq |I| \leq n-d-2$  and $I \subseteq \{d+2, \ldots, n\}$. Among those, the ones that map to $0$ under 
$\pi_{d,n}$ are precisely the $e^1_k$, where $ k=d+2,\dots ,n$ and, moreover, they generate its kernel and their sum is $0$. The images of the rest under 
$\pi_{d,n}$ span the rays of the fan of the projective toric variety $\overline{P}_{d-1,n-1}^{LM}$. 

By the induction hypothesis  $\operatorname{Bl}_{e}\overline{P}_{d-1,n-1}^{LM}$ is not an MDS  and its pseudoeffective cone is not polyhedral,  
so Corollary~\ref{corollary.projection} implies that  $\operatorname{Bl}_{e}\overline{P}_{d,n}^{LM}$ is also not an MDS and its pseudoeffective cone is also not polyhedral.  
The desired conclusion follows by induction.
\end{proof}


We now show that the variety $\operatorname{Bl}_{e}\overline{P}_{d,n}^{LM}$ in Lemma~\ref{lemma.Pdn.Losev.Manin_V2} can be identified with a compactification from Section~\ref{subsection.wonderful}, for an appropriate choice of weights.

\begin{lemma}   \label{lemma.identification}    Consider the set $\mathcal{A'}=(a'_1,\dots,a'_n)\in\mathcal D_{d,n}^{P}$ where $a'_1 = \cdots = a'_{d+1} =1,\, 
a'_{d+2} = \cdots = a'_n =\frac{1}{n-d-2}$. 
Then $Bl_e\overline P_{d,n}^{LM}$ is identified with $\overline P_{d,n}^{\mathcal{A'}}$.
\end{lemma}
\begin{proof}  Let $\mathcal{A}=(a_1=1,\dots,a_{d+1}=1,a_{d+2}=\frac{1}{n-d-1},\dots,a_n=\frac{1}{n-d-1})\in\mathcal D_{d,n}^{P}$ be the set of weights associated to 
$\overline P_{d,n}^{LM}=\overline P_{d,n}^{\mathcal{A}}$.  For $I\subsetneq\{d+1,\dots,n\}$, we have that $H_I\in \mathcal{G}_{\mathcal{A}}$ (resp. $H_I\in\mathcal{G}_{\mathcal{A'}}$) if and only if  $\sum\limits_{i\in I} a_i>1$ (resp. $\sum\limits_{i\in I} a'_i>1$). Therefore, $\mathcal{G}_{\mathcal{A}}=\{ H_I  \subseteq ( \bP^{n-d-2})^d  \; | \;d+1\in I\,\, \text{and}\,\, |I|\geq 2\}$ and $\mathcal{G}_{\mathcal{A'}}=\{ H_I  \subseteq ( \bP^{n-d-2})^d  \; | \;d+1\in I\,\, \text{and}\,\, |I|\geq 2\}\cup\{H_{\{d+2,\dots,n\}}\}=\mathcal{G}_{\mathcal{A}}\cup\{e\}$. Let $\mathcal{G}_{\mathcal{A}}$ be given a total order extending the partial order of ascending dimension. Then, by Proposition~\ref{smoothproj}, we have that $\overline P_{d,n}^{LM}$ is the iterated blowup at $\mathcal{G}_{\mathcal{A}}$ in that order. Since $e$ is disjoint from all $H_I \in \mathcal{G}_{\mathcal{A}}$, the same proposition implies that we may obtain $\overline P_{d,n}^{\mathcal{A'}}$ as a blowup at $\mathcal{G}_{\mathcal{A}}$ followed by blowing up $e$. 
\end{proof}

In the following proposition we show that the blowup $\overline P_{d,n} \rightarrow (\mathbb{P}^{n-d-2})^d$ from Proposition~\ref{smoothproj} factors through $\operatorname{Bl}_{e}\overline P_{d,n}^{LM}$.

\begin{proposition} \label{thm:SecblowUps}
 There exists a sequence of
blowups
\begin{align*}
\overline P_{d,n} \longrightarrow
\operatorname{Bl}_{e}\overline P_{d,n}^{LM}
\longrightarrow
\overline P_{d,n}^{LM}
\longrightarrow
(\mathbb{P}^{n-d-2})^d.
\end{align*}
\end{proposition}
\begin{proof}
By \cite[Proposition~5.1]{GR}, if $\mathcal{A}=(a_1,\dots,a_n)\in\mathcal D_{d,n}^{P}$ 
and $\mathcal{B}=(b_1,\dots,b_n)\in\mathcal D_{d,n}^{P}$
satisfy $a_i \geq b_i$ for all $1 \leq i \leq n$, there exist a
weight reduction morphism $\overline P^{\mathcal A}_{d,n} \rightarrow \overline P^{\mathcal B}_{d,n}$. 
Recall that  we write $\overline P_{d,n}$ in place of $\overline P^{\mathcal A}_{d,n}$ where $\mathcal{A}=(1,1,\dots,1)$. 
The desired conclusion follows from Lemma~\ref{lemma.identification} and the fact that the weight reduction morphisms are given by blowups and are compatible under composition by \cite[Proposition~5.2]{GR}. 
\end{proof}

\vspace{0.1in}

\textit{Proof of Theorem~\ref{thm.cone.pdn}}.
By Proposition~\ref{thm:SecblowUps}, there exists a surjective morphism $\overline{P}_{d,n} \rightarrow \operatorname{Bl}_{e}\overline{P}_{d,n}^{LM}$. 
By Lemma~\ref{lemma.Pdn.Losev.Manin_V2}, $\operatorname{Bl}_{e}\overline{P}_{d,n}^{LM}$ is not a Mori dream space and its pseudoeffective cone is not polyhedral for $n \geq d+9$. Then, the result follows by Theorem~\ref{thm:PropertiesMDS} and Theorem~\ref{thm:PropertiesPseudoeffective}.  
\qed


\subsection{Proof of Theorem~\ref{thm.FM}}        \label{proof.thm.1.1}
 
 The proof of Theorem~\ref{thm.FM} will follow from Theorem~\ref{thm:FM_notMDS} and Proposition~\ref{birFM} below.

\begin{theorem}\label{thm:FM_notMDS}
 $\mathbb{P}^d[n]$  is not a Mori dream space for $n \geq d+9$.
\end{theorem}
\begin{proof}
By Remark~\ref{geomquot},
there exists an $SL_{d+1}$-invariant non-empty open subset of  $\mathbb{P}^d[n]$ such that 
$\overline{P}_{d,n}$ is a good quotient of that open subset by $SL_{d+1}$. 
By Theorem~\ref{thm:bak1}, it follows that  $\mathbb{P}^d[n]$  is not a Mori dream space for $n \geq d+9$ since $\overline{P}_{d,n}$ is not an MDS for $n \geq d+9$ by Theorem~\ref{thm.cone.pdn}.
\end{proof}

\begin{proposition}\label{birFM} Let $f:X\dashrightarrow Y$ be a birational map between two smooth projective varieties $X$ and $Y$. Suppose that there exists an open subset $U$ of $X$ and an open subset $V$ of $Y$
with $codim(Y\setminus V,Y)\geq 2$ such that $f_{|U}$ is an isomorphism onto $V$. If $X[n]$ is a Mori dream space for some $n>0$, then $Y[n]$ is a Mori dream space.
\end{proposition}
\begin{proof} By the assumptions, $V^n$ is an open subset whose complement has codimension $\geq 2$ in $Y^n$. Let $I\subseteq \{1,\dots,n\}$ such that $|I|\geq 2$. Observe that the intersection of $V^n$ with the diagonal in $Y^n$ corresponding to $I$ is the diagonal in $V^n$ corresponding to $I$. In other words, each center in the iterated blowup construction of $Y[n]\to Y^n$ restricts to a center in the iterated blowup construction of $V[n]\to V^n$ under the inclusion $V^n\to Y^n$ and every center of $V[n]\to V^n$ is obtained via this restriction. Since blowing up commutes with open immersions, we obtain an open immersion of the Fulton-MacPherson space $V[n]$ in $Y[n]$ whose complement in $Y[n]$ has codimension at least two. Similarly, the Fulton-MacPherson space $U[n]$ is an open subvariety of $X[n]$. By the hypothesis, it follows that $U[n]$ is isomorphic to $V[n]$. 

To conclude the claim, we argue by contradiction: let us suppose that $Y[n]$ is not an MDS. This implies that its Cox ring is not finitely generated, and since $codim(Y[n]\setminus V[n],Y[n])\geq 2$, the Cox ring of $V[n]$ is not finitely generated. Since $U[n]$ is isomorphic to $V[n]$, we deduce that the Cox ring of $U[n]$ is not finitely generated either. 
Let $\Gamma$ be a finitely generated  semigroup of Weil divisors on $U[n]$ such that their classes span the vector space  $Cl(U[n])_{\mathbb{Q}}$. 
Let us denote by $\overline{\Gamma}$ the image of $\Gamma$ under the $\mathbb{Z}$-linear inclusion of $WDiv(U[n])$ into $WDiv(X[n])$ induced by sending a prime divisor of $U[n]$ to its closure in $X[n]$. 
Since $X[n]$ is an MDS, the algebra 
\[
\bigoplus_{D \in \overline{\Gamma}}H^{0}(U[n],\mathcal{O}_{X[n]}(D)) \cong \bigoplus_{D \in \Gamma}H^{0}(U[n],\mathcal{O}_{U[n]}(D))  
\]
is finitely generated by Theorem~\ref{thm:bak1B}.
However, this algebra is a Cox ring of $U[n]$, so it cannot be finitely generated, which gives us a contradiction. 
\end{proof}

As an immediate consequence of the above proposition
and Theorem~\ref{thm:FM_notMDS}, we get the following.
\begin{corollary} Let $X$ be a smooth projective variety. Suppose that there exist smooth projective varieties
$X_0:=\PP^d,X_1,\dots, X_k:=X$ for some $k>0$ and a sequence of birational maps $$X= X_k\dashrightarrow X_{k-1}\dashrightarrow\dots \dashrightarrow X_0=\PP^d$$
where $X_i\dashrightarrow X_{i-1}$ is either a blowup morphism or a small modification for $i=1,\dots ,k$. 
Then $X[n]$ is not a Mori dream space for $n\geq d+9$.
\end{corollary}

\bibliographystyle{amsalpha}
\bibliography{FultonMacPhersonNotMDS}

\end{document}